\documentclass[12pt,a4paper]{amsart}
\newtheorem{de}{Definition}[section]

\newtheorem{re}[de]{Remark}

\newtheorem{te}[de]{Theorem}

\usepackage{soul}
\usepackage{amsfonts} 
\usepackage{amsmath}
\usepackage[all]{xy}

\newcommand{\ot}{\otg{B}}

\newcommand{\mproof}{\noindent{\bf Proof.}}
\newcommand{\twosid}[3]{\ar@<0.25ex>@{<-}[#1]^{#2}  \ar@<-1ex>[#1]_{#3}}

\newenvironment{Proof}{{\mproof}}{{ }\hfill$\Box$ \newline\vspace{\lineskip}}

\newcommand{\sw}[1]{{}_{(#1)}}

\def\eps{\varepsilon}

\def\id{\textrm{Id}}

\def\ot{\otimes}

\begin{document}

\title{On Jordan Algebras and Some Unification Results}
\maketitle
\begin{center}
\author{Florin F. Nichita}\\
{\it Simion Stoilow Institute of Mathematics
of the Romanian Academy\\
21 Calea Grivitei, 010702 Bucharest, Romania}
\end{center}

\bigskip

\bigskip

\begin{center}
 {\bf Abstract}
\end{center}

This paper is based on a talk given at the 14-th International Workshop on Differential
Geometry and Its Applications, 
hosted by the Petroleum Gas University from Ploiesti, between
July 9-th and July 11-th, 2019.
After presenting some historical facts,
we will consider some geometry problems
related to unification approaches.
Jordan algebras and Lie algebras are the main non-associative structures.
Attempts to unify non-associative algebras
and associative algebras
led to UJLA structures. 
Another algebraic structure which unifies
non-associative algebras
and associative algebras
is the Yang-Baxter equation.
We will review topics relared to the
Yang-Baxter equation and Yang-Baxter systems,
with the goal to unify constructions
from Differential Geometry.



{ \bf MSC-class:} 01A99, 00A05, 16T25 (Primary), 17C05, 17C50 (Secondary)


\section{Introduction}

This paper is a survey paper, but it 
 also presents new results. It
emerged after our lectures given at the 9--th Congress of 
Romanian Mathematicians (CRM 9, Galati, Romania, June-July 2019) and
14-th International Workshop on Differential
Geometry and Its Applications
(DGA 14, UPG Ploiesti, July 2019).

We begin with some historical remarks.
Vr\u anceanu proposed the study of spaces 
with constant affine connection
associated to finite-dimensional real Jordan algebras in 1966.
This study was 
afterwards developed by Iordanescu, 
Popovici and Turtoi (see the Comments on page 31 from \cite{RI}).
Dan Barbilian 
 proved 
that the rings which can be the underlying rings for projective geometries are
(with a few exceptions) rings with a unit element in which any one-sided inverse
is a two-sided inverse  (see \cite{RI}).
Another big Romanian geometer was G. Tzitzeica, and 
some of his contributions to mathematics
will be 
presented in our next section.

A substantial part of this paper is dedicated to unifying structures.
The apparition of the Yang-Baxter equation (\cite{py}) in theoretical physics and statistical mechanics
(see \cite{yang, baxter1, baxter2}) has led to many applications in these fields and in 
quantum groups, quantum computing, knot theory, braided categories, analysis of integrable systems, quantum mechanics, etc (see \cite{si}).
The interest in this equation is growing, as new properties of it are found, and
its solutions are not classified yet.
This equation can be understood as an 
unifying equation (see, for example, \cite{fnichita, lebed, lebed2, inn2}).
Another unification of non-associative structures, was 
obtained using the UJLA structures
(\cite{inn2, inn}),
 which 
could be seen as structures which comprise the information
encapsulated in associative algebras, Lie algebras and Jordan algebras.

We now refer to important results obtained under the 
guidance of St. Papadima. Some presentations at DGA 14 were dedicated
to his memory (see \cite{Paolo, RaduP, Florin, Elis}). It was observed
during DGA 14, that some results from \cite{RaduP} can be extended for the
virtual braid groups from \cite{Paolo}. Also, the techniques of \cite{Elis}
can be considered in the framework of UJLA structures.
The Yang--Baxter equation  plays an important role in knot theory. 
D\u asc\u alescu and Nichita have shown in \cite{DN} how to associate a Yang--Baxter operator
to any algebra structure over a vector space, using the associativity 
of the multiplication. 
 Turaev has described in \cite{T} 
a general scheme to derive an invariant of oriented links 
from a Yang--Baxter operator, provided this one can be ``enhanced''.
The invariant which was obtained from
those Yang--Baxter operators is the Alexander polynomial of knots.
Thus, in a way, the Alexander polynomial is the knot invariant corresponding to
the axioms of (unitary associative) algebras. 

\bigskip
The current paper is organised as follows. 
The next section deals with  
Tzitzeica-Johnson's theorem.
In Sections 3 and 4 we present the Yang-Baxter
equations and  the Yang-Baxter systems(including the
set-theoretical Yang-Baxter
equation). 
These will be followed by a section on Jordan algebras, 
UJLA structures and unification constructions in Differential
Geometry.
A conclusion section will end this article.  

\bigskip

\bigskip

{\centering \section{On Tzitzeica-Johnson's theorem and pictorial mathematics}}

In his talk at DGA 12,
Florin Caragiu explained that there 
exists a special mathematical discourse, called
``proofs without words'', which uses pictures or diagrams 
in order to boost the intuition of the reader 
(see \cite{Cara}). Pictorial (diagrammatic) style 
of mathematical language is much appreciated by both educators and 
researchers in mathematics (see, also, \cite{Tzitzeica}). 
Very easy to be grasped, some pictorial style problems 
need an entirely ``artillery'' in order to be cracked.

\bigskip

We now consider the  Tzitzeica-Johnson's problem (see, for example,
 \cite{Tzitzeica}). 
We start with three circles of the same radius $ r$.
 The intersection points
of pairs of circles 
are denoted by A, B, C and the common point of intersection of the three circles is 
denoted by O (see the Figure 1 below).

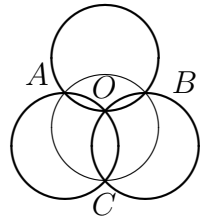
\begin{figure}[h!]
  \caption{{\bf The three coins problem} (Tzitzeica - 1908, Johnson - 1916). Under 
the above assumptions, there exists a circle with radius $r$, passing through
the points A, B and C.}
\begin{picture}(140,87)
\put(145,37){\circle{50}}
\thicklines
\put(145,63){\circle{50}}
\put(160,30){\circle{50}}
\put(130,30){\circle{50}}
\put(115,53){$A$}
\put(140,48){$O$}
\put(170,50){$B$}
\put(140,4){$C$}
\end{picture}
\end{figure}

We now propose a new construction (in the
Figure 2).

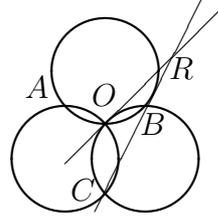
\begin{figure}[h!]
  \caption{{The tangent to the third circle in O meets the line BC in R.}}
\begin{picture}(140,90)
\put(130,28){\line(1,1){60}}
\put(141,10){\line(1,2){40}}
\thicklines
\put(145,63){\circle{50}}
\put(160,30){\circle{50}}
\put(130,30){\circle{50}}
\put(115,53){$A$}
\put(140,50){$O$}
\put(169,60){$R$}
\put(158,39){$B$}
\put(132,14){$C$}
\end{picture}
\end{figure}

\bigskip
\bigskip
\bigskip

At this moment, we can present our {\bf Theorem}: under
the above assumptions, 
the tangents to the circles in O meet the lines
AB, AC and BC in three colinear points (see the
Figure 3 below).

\begin{figure}[h!]
  \caption{{{\bf Theorem.} 
The points P, Q and R are colinear.}}
\begin{picture}(140,100)
\put(110,49){\line(1,0){90}}
\put(130,28){\line(1,1){60}}
\put(141,10){\line(2,5){30}}
\put(130,39){\line(5,1){60}}
\put(180,8){\line(-1,1){84}}
\put(148,10){\line(-3,5){55}}
\thicklines
\put(140,62){\circle{50}}
\put(160,30){\circle{50}}
\put(130,30){\circle{50}}
\put(107,81){$Q$}
\put(170,50){$P$}
\put(149,59){$R$}
\end{picture}
\end{figure}
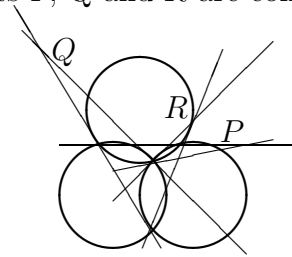

\bigskip

The proof of the above theorem is based on the properties of the power of a circle
and on the Desargues' theorem. It is
 the limit case of the following
picture. 

\begin{figure}[h!]
  \caption{{The limit case of the following picture
implies our theorem.}}
\begin{picture}(140,86)
\put(145,10){\line(0,1){65}}
\put(125,16){\line(6,5){55}}
\put(165,16){\line(-6,5){55}}
\thicklines
\put(145,48){\circle{50}}
\put(160,30){\circle{50}}
\put(130,30){\circle{50}}
\end{picture}
\end{figure}
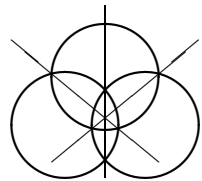

The above problems can be included in a more general scheme 
(see the Figure 5).

\begin{figure}[h!]
  \caption{{\bf A mathematical scheme.} If $\mathcal{X}$, $\mathcal{Y}$ and
$\mathcal{Z}$ have certain properties, then there exists a related
$\mathcal{T}$, with the same properties.}
\begin{picture}(140,86)
\thicklines
\put(145,48){\circle{50}}
\put(160,30){\circle{50}}
\put(130,30){\circle{50}}
\put(115,65){$\mathcal{X}$}
\put(182,30){$\mathcal{Z}$}
\put(99,24){$\mathcal{Y}$}
\end{picture}
\end{figure}
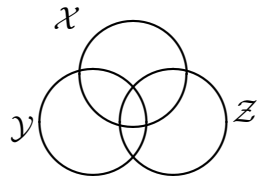

\bigskip
\bigskip

\newpage

The dual of the Desargues' Theorem can be interpretated in the light of the above scheme by
using arrangements of (three) colored lines.

The hypothesis:

$\mathcal{X} = \{ (d, f, h) \ such \ that \ d \cap f \cap h = P \}$

$\mathcal{Y} = \{ (d', f, h') \ such \ that \ d' \cap f \cap h' = Q \}$

$\mathcal{Z} = \{ (d'', f, h'') \ such \ that \ d'' \  \cap f \cap h'' = R \}$

$\mathcal{X} \cap \mathcal{Y} \cap \mathcal{Z} = \{ f \} $

$\mathcal{X} \cap \mathcal{Y} = \{ A, A' \} \cup \{ f \} $

$\mathcal{Y} \cap \mathcal{Z} = \{ B, B' \}  \cup \{ f \}  $

$\mathcal{X} \cap \mathcal{Z} = \{ C, C' \} \cup \{ f \} $

Conclusion:

There exist a triple $ ( l, m, n )$ through $ A, A', B, B', C, C'$.

\bigskip

From the phylosophical point of  view,
one could include in this scheme the 
{ \it inclussion-exclusion principle}:

$ \ \vert {X} \cup {Y} \cup  {Z} \vert =
\vert {X} \vert + \vert {Y} \vert +
\vert {Z} \vert - \vert {X} \cap {Y} \vert -
\vert {X} \cap {Z} \vert -
\vert {Y} \cap {Z} \vert + \vert {X} 
\cap {Y}  \cap {Z} \vert $.

Let us explain this analogy.
Once we have information about the cardinality of 
$ X, \ Y, \ Z, \ {X} \cap {Y} , \ {X} \cap {Z}, \ {Y} \cap {Z}
$ and $ {X} 
\cap {Y}  \cap {Z} $,
then there exists a $T$, related to
$ X, \ Y$ and $  Z$,
 whose cardinality can be also computed. Of course, in this case
$T= {X} \cup {Y} \cup  {Z} $.

In order to illustrate better this analogy, let us
suppose further that 
$\vert {X} \vert = \vert {Y} \vert =
\vert {Z} \vert = 7$,
$\vert {X} \cap {Y} \vert =
\vert {X} \cap {Z} \vert =
\vert {Y} \cap {Z} \vert= 1$
and
$\vert {X} 
\cap {Y}  \cap {Z} \vert = 1 $.
Then, there exists a set $T$, with $\vert {T} \vert = 7 $
which contains elements from all possible (sub)sets.
\bigskip

Notice that Tzitzeica-Johnson's Theorem has many generalizations 
(\cite{Tzitzeica}) and interpretations
related to the above mathematical scheme (from Figure 5).

 Thus, Tzitzeica-Johnson's Theorem can  be
interpretated in terms of disks: 

(i) ``If three disks
of the same radius have a common point of intersection, then they contain inside of
their union a forth disk with the same radius.'' 

(ii) ``If three disks $ \mathcal{X} $, $ \mathcal{Y} $ and $ \mathcal{Z} $
of the same radius have a common point, then there exists a forth disk
of the same radius which includes  $( \mathcal{X} \cap \mathcal{Y}) \cup  
(\mathcal{X} \cap \mathcal{Z}) \cup  
(\mathcal{Y} \cap \mathcal{Z})$.''

It is an open problem to prove
a similar statement for a domain bounded by
an arbitary closed convex curve.
In other words we conjecture that if we consider 
a domain $ \mathcal{X} $
 bounded by
an arbitary closed convex curve,
and two copies of  $ \mathcal{X} $, denoted by
 $ \mathcal{Y} $ and  $ \mathcal{Z} $ such that
$ \mathcal{X} \cap 
\mathcal{Y}
 \cap 
\mathcal{Z} \neq \emptyset $, then there exists another copy of
 $ \mathcal{X} $, denoted by
 $ \mathcal{T} $ such that
 $ \ \ ( \mathcal{X} \cap \mathcal{Y}) \cup  
(\mathcal{X} \cap \mathcal{Z}) \cup  
(\mathcal{Y} \cap \mathcal{Z}) 
  \subset \mathcal{T}
$.
\bigskip
  
\bigskip
We now consider another type of problems.

In his best-seller book  \cite{Greco}, P. Greco (see \cite{pi} - there
is a {\it puzzle} there!) refers
to a problem which can be generalized as follows.

\bigskip

{\it OPEN PROBLEMS. For an arbitrary convex closed curve, 
we consider the length of largest diameter (D). 
Here, the largest diameter
 is the longest
segment through the center of mass 
with the endpoints on the given curve.
 In a similar manner, one can define the smallest 
diameter (d). 

i) If $L$ is the length of the given curve
and the domain inside the
given curve is a convex set, then
we conjecture that:
$$ \frac{L}{D} \leq \pi \leq \frac{L}{d} \ \ .$$

(ii) If $A$ is the area inside the given curve, the equation 

\begin{equation} \label{ecgraddoi}
x^2 - \frac{L}{2} x + A =0 
\end{equation}
and its implications are not completely understood. For example,
if the given curve is an ellipse, solving this equation 
in terms of the semi-axes of the ellipse is an unsolved problem. }

\bigskip
\bigskip

The second part of the above conjecture is related to \cite{nicu}.

\bigskip

Another idea related to \cite{Greco} and to
the Euler's relation,
  $ \ e^{ \pi i} + 1 = 0 \ $, is the following
 relation containing $ \pi$:
$ \ \vert e^i - \pi \vert >  e $.
It is very good approximation, and it has a pictorial interpretation
(see the Figure 6).

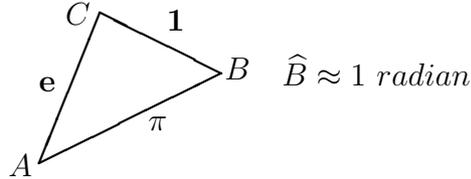
\begin{figure}[h!]
  \caption{{An interpretation for the formula 
$ \  \vert e^i - \pi \vert \approx  e \ $.}}
  \centering
\setlength{\unitlength}{0.8cm}
\begin{picture}(4,3.2)
\thicklines
\put(1,0.5){\line(2,1){3}}
\put(4,2){\line(-2,1){2}}
\put(2,3){\line(-2,-5){1}}
\put(0.5,0.3){$A$}
\put(4.05,1.9){$B$}
\put(1.44,2.80){$C$}
\put(3.1, 2.7){${\bf 1}$}
\put(1, 1.7){${\bf e}$}
\put(2.8, 1.05){${\bf \pi}$}
\put(5,1.8){$\displaystyle
\widehat{B} \approx 1 \ radian$}
\end{picture}
\end{figure}

\section{Yang-Baxter equations}

In the next sections, we will work over a generic field $ k $.
The tensor products will be defined over $k$.

This section was inspired by two papers,
\cite{mnp, SOKT}, and it
begins with an introduction to the Yang-Baxter equation.

\bigskip

Let $ V $ be a vector space over $ k$.
Let $ I=I_V : V \rightarrow V $ be
the identity map of the space $V$.

 We denote by
$ \   \tau : V \otimes V \rightarrow V \ot V \  $ the twist map defined by $ \tau (v \ot w) = w \ot v $.

\bigskip

Now, we will give the main notations for introducing the Yang-Baxter equation.

For $ \  R: V \ot V \rightarrow V \ot V  $
a $ k$-linear map, let
$ {R^{12}}= R \ot I , \  {R^{23}}= I \ot R : V \ot V \ot V \rightarrow V \ot V\ot V \ .
$
In a similar manner, we denote by
${R^{13}}$ a linear map acting on the first and third component of
$ V \otimes V \otimes V$. It turns out that
${R^{13}}=(I\ot\tau )(R\ot I)(I\ot \tau ) $.

\bigskip

\begin{de} A
{ \bf  Yang-Baxter
operator} is $ k$-linear map
$ R : V \ot V \rightarrow V \ot V $,
which satisfies the braid condition (the Yang-Baxter equation):
\begin{equation}  \label{ybeq}
R^{12}  \circ  R^{23}  \circ  R^{12} = R^{23}  \circ  R^{12}  \circ  R^{23}.
\end{equation}
We also require that the map $R$ is invertible.
\end{de}

\begin{re}
 Some examples of  Yang-Baxter
operators are the following:
$ R = \tau \ $ (i.e., $ \ R(a \otimes b) = b \otimes a$), and
$ R = I \otimes I \ $ (i.e., $ \ R(a \otimes b) = a \otimes b$).
\end{re}

\begin{re}
An important observation is that if $R$ satisfies (\ref{ybeq}) then both
$R\circ \tau  $ and $ \tau \circ R $ satisfy the QYBE
(the quantum Yang-Baxter equation):
\begin{equation}   \label{ybeq2}
R^{12}  \circ  R^{13}  \circ  R^{23} = R^{23}  \circ  R^{13}  \circ  R^{12}.
\end{equation}
\end{re}

\begin{re}
The equations (\ref{ybeq}) and (\ref{ybeq2}) are equivalent.
\end{re}

\begin{re}

The following construction of Yang-Baxter operators is described in \cite{DN}.

If $A$ is a $k$-algebra, then for all non-zero $r,s\in k$, the linear map 
\begin{equation}\label{ra}
R^A_{r,s}:A \ot A
\rightarrow A \ot A,\quad   a \ot b \mapsto s ab \ot 1 + r 1 \ot ab - s a \ot b 
\end{equation}
is a Yang-Baxter operator. 

The inverse of $R^A_{r,s}$ is 
${(R^A_{r,s})}^{-1}(a \ot b)= \frac{1}{r} ab \ot 1 + \frac{1}{s} 1 \ot 
ab - \frac{1}{s} a \ot b $.

\end{re}

 Turaev has described in \cite{T} 
a general scheme to derive an invariant of oriented links 
from a Yang--Baxter operator, provided this one can be ``enhanced''.

The Jones polynomial \cite{Jones} and its two--variable extensions, 
namely the Homflypt polynomial \cite{HOMFLY,PT} and the Kauffman polynomial \cite{K}, 
can be obtained in that way by ``enhancing'' some Yang--Baxter operators 
obtained in \cite{Jimbo}. Those solutions of the Yang--Baxter equation are associated to simple Lie algebras 
and their fundamental representations. The Alexander polynomial can be derived
from a Yang--Baxter operator as well, using a slight 
modification of Turaev's construction \cite{M}.

The only invariant which can be obtained from the
Yang--Baxter operators 
(\ref{ra})
 is the Alexander polynomial of knots.
Thus, in a way, the Alexander polynomial is the knot invariant corresponding to
the axioms of (unitary associative) algebras (cf. \cite{mn}). 

Note that specializations of the Homflypt polynomial had to be expected 
from those Yang--Baxter operators since they have degree $2$  minimal polynomials.

\bigskip

\begin{re}
There is a similar terminology for the set-theoretical Yang-Baxter equation. In this
case $V$ is replaced by a set $X$ and the tensor product by the Cartesian product.
We will explain this definition in the next examples below.
\end{re}


Let $ X $  be a set containing
three logical sentences $p, \ q, \ r$ (i.e., $ p, \ q, \ r  \in X$).
We can choose $ X $ as rich as we wish for the moment. Later, we will
try to find the smallest $ X $ which fits for our theory.


Let $R  : X \times X \rightarrow X \times X $, be defined by
$ R( p, \ q) = ( p \vee q, \  p \wedge q) $.

It follows that
\begin{equation}  \label{stybeq}
  (R \times I)  \circ (I \times R)  \circ  (R \times I) =
(I \times R)  \circ  (R \times I)  \circ (I \times R) \ .
\end{equation}

One way to check that  (\ref{stybeq}) holds is to make a table with values
for $p, \ q, \ r$.

\bigskip

\begin{re}

Another interesting solution to the  set-theoretical Yang-Baxter equation
is the following.

Let $R  : X \times X \rightarrow X \times X $, be defined by
$ R( p, \ q) = ( p \rightarrow q, \  p ) $.

Again, one way to check the above statement is to make a table with values
for $p, \ q, \ r$.
\end{re}

\bigskip

Let us denote $ R( p, \ q) = ( p_1, \  q_2) $.
The set-theoretical Yang-Baxter equation can be expressed as

\begin{equation}\label{unfa}
 ( p_{101}, \ q_{212}, \ r_{020}) =
(p_{010}, \ q_{121}, \ r_{212} ) \ .
\end{equation}

\bigskip
It is now the moment to discuss about  $X$.

{\bf Open problems.}
What can be said about $X$ in general?
What is the smallest $X$ for which $R$ is well-defined?
One could consider a set $X$ containing more than three logical sentences.
What is the interpretation of the set-theoretical
Yang-Baxter solutions in this case ?

We can go a step further and consider an algebra of type (2,2),
$ (A, *, \circ )$, and call the operations $*$ and $\circ$
 {\it YB conjugated}
if $ R(a,b)= ( a*b, a \circ b) $ satisfies the set-theoretical Yang-Baxter
equation.

We propose the study of algebras with two operations  which
are YB conjugated. 
Groups,
distributive lattices,  
and self-distributive laws can be considered 
as objects with operations which are YB conjugated
(compare with \cite{tom}).

\bigskip

\begin{re} We now consider other (systems of) equations for a $ k$-linear map
$ R : V \ot V \rightarrow V \ot V $:
\begin{equation}  \label{i}
R^{12}  \circ  R^{23}  = R^{13}  \circ  R^{12} =   R^{23} \circ R^{13}
\end{equation}
\begin{equation}  \label{x}
R^{23}  \circ  R^{12}  = R^{12}  \circ  R^{13} 
\end{equation}
\begin{equation}   \label{nou}
R^{12}  \circ  R^{13}  \circ  R^{12} \circ  R^{23}
= R^{13}  \circ  R^{23}  \circ  R^{13} \circ  R^{12}
\end{equation}
\begin{equation}  \label{unifyb}
(R^{12}  \circ  R^{23}  \circ  R^{12} - 
R^{23}  \circ  R^{12}  \circ  R^{23}) \circ
(R^{12}  \circ  R^{13}  \circ  R^{23} -
 R^{23}  \circ  R^{13}  \circ  R^{12}) =0
\end{equation}
\begin{equation}   \label{nou2}
R=XY \ ,  \ \
R^{12}  \circ  X^{13}  \circ  R^{23} \circ  Y^{12}
= R^{23}  \circ  X^{13}  \circ  R^{12} \circ  Y^{23}
\end{equation}
\begin{equation}   \label{2id}
R \circ R = I \otimes I
\end{equation}
\begin{equation}  \label{unifyb2}
R^{12}  R^{23}      R^{13}    R^{23} +
R^{23}    R^{12}      R^{13}    R^{12} = 
R^{23}    R^{12}   R^{23} 
R^{12}    R^{13}    R^{23} +
R^{12}    R^{23}    R^{12} 
R^{23}    R^{13}    R^{12}
\end{equation}

\end{re}

\begin{te}
 If a $ k$-map
$ R : V \ot V \rightarrow V \ot V $ verifies (\ref{i}) and (\ref{x})
then $R$ is a common solution for (\ref{ybeq}) and (\ref{ybeq2}).

If a $ k$-map
$ R : V \ot V \rightarrow V \ot V $ verifies
(\ref{ybeq}) and (\ref{ybeq2}) then it is a solution for
 (\ref{nou}).

If a $ k$-map
$ R : V \ot V \rightarrow V \ot V $ verifies
(\ref{ybeq}) or (\ref{ybeq2}) then $R$ is a solution for
 (\ref{unifyb}).

If a $ k$-map
$ R : V \ot V \rightarrow V \ot V $ verifies
(\ref{ybeq}) or (\ref{ybeq2}) then it is a solution for
 (\ref{nou2}).

\end{te}

{\bf Proof.} We only prove the first claim:
$ R^{23}  \circ  R^{12}  \circ  R^{23} =
R^{23}  \circ  R^{13}  \circ  R^{12} = 
 R^{12}  \circ  R^{23}  \circ  R^{12} $ 

$
 R^{23}  \circ  R^{13}  \circ  R^{12}=
R^{23}  \circ  R^{12}  \circ  R^{23}=
R^{12}  \circ  R^{13}  \circ  R^{23} $

The other claims follow in a similar manner. $\diamond$
\bigskip
\bigskip

We will write the above results in the following manner: 
\begin{center}
$ (\ref{i}) \wedge (\ref{x}) \rightarrow 
(\ref{ybeq}) \wedge (\ref{ybeq2})$;
$ \ \ \ \ \ (\ref{ybeq}) \wedge (\ref{ybeq2}) \rightarrow 
(\ref{nou})$;

$ (\ref{ybeq}) \vee (\ref{ybeq2}) \rightarrow 
(\ref{unifyb})$;
$ \ \ \ \ \ (\ref{ybeq}) \vee (\ref{ybeq2}) \rightarrow 
(\ref{nou2})$; $ \ \ \ (\ref{unifyb}) \wedge (\ref{2id})
\rightarrow (\ref{unifyb2})$.

\end{center}

As a direct application of this theorem, one can check which of the
funtions presented in this section are common solutions for the
braid condition and the quantum Yang-Baxter equation.

\bigskip

\begin{re}
 The set of equations of operators of type
$ R : V \ot V \rightarrow V \ot V $ and
$ R^{ij} : V \ot V \ot V
\rightarrow V \ot V \ot V $, with $ \ i,j \in \{1, 2, 3 \} $,
has a natural distributive latice structure (see \cite{mnp}).

\end{re}


\section{ Yang-Baxter systems}

We present the Yang-Baxter systems theory following the paper \cite{bn}.
Yang-Baxter systems were introduced in \cite{HlaSno:sol} as a 
spectral-parameter independent generalisation of 
quantum Yang-Baxter equations related to non-ultralocal 
integrable systems. 

Yang-Baxter systems are conveniently defined in 
terms of {\em Yang-Baxter commutators}. 
Consider three vector spaces $V,V',V''$ and  three linear maps
$ R : V \ot V' \rightarrow V \ot V' $,
$ S : V \ot V'' \rightarrow V \ot V'' $ and
$ T : V' \ot V'' \rightarrow V' \ot V'' $. 

Then a {\em Yang-Baxter commutator} is a map
$ [R,S,T]:  V \ot V' \ot V'' \rightarrow V \ot V' \ot V'' $, defined by
\begin{equation}   \label{ybcomm}
[R,S,T]= R_{12}  \circ  S_{13}  \circ  T_{23} - T_{23}  \circ  S_{13} 
\circ  R_{12} \ .
\end{equation}

In terms of a Yang-Baxter commutator,  the quantum Yang-Baxter equation (\ref{ybeq2})  is expressed simply as $[R,R,R] = 0$.

\begin{de}\label{def.wxz}
Let $V$ and  $V'$ be vector spaces. A system of linear maps 
$$ W : V \ot V \rightarrow V \ot V ,\quad 
Z : V' \ot V' \rightarrow V' \ot V' , \quad
 X : V \ot V' \rightarrow V \ot V' 
$$ 
is called  a {\em WXZ-system} or a 
{\em Yang-Baxter system}, provided the following equations are satisified:
\begin{equation}   \label{ybeqn4}
[W,W,W]\ = \  0 \ ,
\end{equation}
\begin{equation}   \label{ybeqn5}
[Z,Z,Z]\ = \  0 \ ,
\end{equation}
\begin{equation}   \label{ybeqn6}
[W,X,X]\ = \  0 \ ,
\end{equation}
\begin{equation}   \label{ybeqn7}
[X,X,Z]\ = \  0 \ .
\end{equation}
\end{de}

\bigskip

There are several algebraic origins and applications of WXZ-systems. 
 WXZ-systems with invertible $W$, $X$ and $Z$ can be used to 
construct dually-paired bialgebras of the FRT type, thus leading to 
quantum doubles. 

More precisely, consider a WXZ-system with finite-dimensional $V=V'$, so that each of $W$, $X$, $Z$ is an $N^2\times N^2$-matrix. Suppose that $W,X,Z$ are invertible. Since $W$ and $Z$ satisfy Yang-Baxter equations (\ref{ybeqn4})--(\ref{ybeqn5}), one can consider two matrix bialgebras $A$ and $B$ with $N\times N$ matrices of generators $U$ and $T$ respectively, and relations
$ W_{12} U_{1} U_{2} = U_{2} U_{1} W_{12} $, 
$ Z_{12} T_{1} T_{2} = T_{2} T_{1} Z_{12}$. 
The existence of 
an invertible operator $X$ that satisfies equations (\ref{ybeqn6})--(\ref{ybeqn7}), 
means that $A$ and $B$ are dually paired with a non-degenerate pairing $ <U_{1},T_{2}> = X_{12} $. Furthermore,  the tensor product $A\ot B$ has an algebra (quantum double) structure with crossed relations
$ X_{12}U_{1}T_{2} = T_{2}U_{1}X_{12} $.

Given a WXZ-system as in Definition~\ref{def.wxz} one can 
construct a Yang-Baxter operator on $V\oplus V'$, provided the map $X$ 
is invertible. This is a special case of a {\em gluing procedure} 
described in  \cite[Theorem~2.7]{MajMar:glu} 
(cf.\  \cite[Example 2.11]{MajMar:glu}).
Let $ R = W \circ \tau_{V,V} $, $ R' = Z \circ \tau_{V',V'} $, $ U = X
\circ \tau_{V',V} $. Then the linear map
$$
 R \oplus_{U} R': (V \oplus V') \ot (V \oplus V') \to  (V \oplus V') \ot (V \oplus V')
$$ 
given by
$ R \oplus_{U} R'|_{V \ot V} = R $,
$ R \oplus_{U} R'|_{V' \ot V'} = R' $, and for all $ x \in V $, $ y \in V' $,
$$ 
(R \oplus_{U} R')(y \ot x) = U(y \ot x ), \qquad  (R \oplus_{U} R')(x \ot y ) = U^{-1}(x \ot y)
$$
is a Yang-Baxter operator.

Entwining structures were introduced  in order 
to recapture the symmetry structure of 
non-commutative (coalgebra) principal bundles or 
coalgebra-Galois extensions. 


\begin{de}\label{def.entw} 
 An algebra $A$ is said to be {\em entwined} with  a coalgebra $C$ if there exists  a linear  map 
$ \psi : C \ot A
\rightarrow A \ot C $ satisfying the following four conditions:

(1) $ \psi \circ (I_{C} \ot \mu) = ( \mu \ot I_{C}) \circ (I_{A} \ot
\psi) \circ ( \psi \ot I_{A}) $,

(2) $ (I_{A} \ot \Delta) \circ \psi = ( \psi \ot I_{C} ) \circ (I_{C}
\ot \psi) \circ ( \Delta \ot I_{A} ) $,

(3) $ \psi \circ ( I_{C} \ot \iota ) = \iota \ot I_{C} $,

(4) $ (I_{A} \ot \eps ) \circ \psi = \eps \ot I_{A} $.

The map $ \psi $ is known as an {\em entwining map}, and the triple  $ {(A,C)}_{\psi} $ is called an {\em entwining structure}. 
\end{de}

\begin{re}

To denote the action of an entwining
map $ \psi $ on elements it is convenient to use 
the following {\em $\alpha$-notation}, for all $a,b\in A$ and $c\in C$,

$
\psi (c \ot a) = \sum_{\alpha} a_{\alpha} \ot c^{\alpha} , \quad
(I_{A} \ot \psi) \circ ( \psi \ot I_{A}) (c \ot a \ot b) =
\sum_{\alpha, \beta} a_{\alpha} \ot b_{\beta} \ot c^{\alpha \beta},
$
etc. 

For example (\ref{unfa}) is a kind of {\em $\alpha$-notation} related to
the Yang-Baxter equation.

The relations (1), (2), (3)
and (4) in Definition~\ref{def.entw} are  equivalent to the following explicit relations, for all $ a,b
\in A$, $c \in C$,
\begin{equation}\label{a}
\sum_{\alpha}(ab)_{\alpha}\ot c^{\alpha} = \sum_{\alpha, \beta}
a_{\alpha} b_{\beta} \ot c^{\alpha \beta},
\end{equation}
\begin{equation}\label{b}
\sum_{\alpha} a_{\alpha} \ot {c^{\alpha}}\sw{1} \ot {c^{\alpha}}\sw{2} =
\sum_{\alpha, \beta} a_{\beta \alpha} \ot {c\sw{1}}^{\alpha} \ot
{c\sw{2}}^{\beta},
\end{equation}
\begin{equation}\label{c}
\sum_{\alpha} 1_{\alpha} \ot c^{\alpha} = 1 \ot c,
\end{equation}
\begin{equation}\label{d}
\sum_{\alpha}a_{\alpha} \eps (c^{\alpha}) = a \eps (c).
\end{equation}
\end{re}

\begin{te}\label{thm.main} (\cite{bn})
 Let $A$ be an algebra  and let $C$ be a
coalgebra. For any $ s, r, t, p \in k$ define linear maps 
$$W: A\ot A\to A\ot A, \qquad a\ot b \mapsto s ba \ot 1 + r 1 \ot ba - s b \ot a, 
$$
$$
Z: C\ot C\to C\ot C, \qquad c\ot d \mapsto t \eps (c) \sum d\sw 1 \ot d\sw 2 + p \eps (d) 
\sum c\sw 1 \ot c\sw 2 - p d \ot c.
$$
Let $X:A\ot C\to A\ot C$ be a linear map such that $ X\circ (\iota \ot \id_C) = \iota\ot \id_{C} $ and $ (\id_A \ot \eps ) \circ X =\id_A\ot \eps$. Then $W,X,Z$ is a Yang-Baxter system if and only if $A$ is entwined with $C$ by the map $\psi: = X\circ \tau_{C, A}$.

\end{te}

\bigskip

{\centering\section{Jordan algebras and UJLA structures}}

Jordan algebras have applications
in physics, differential geometry, ring geometries, quantum groups,
analysis, biology, etc (see \cite{I, RI, RIordanescu}).

We have introduced some structures which unify Jordan algebras, Lie algebras
and (non-unital) associative algebras. 
These structures were called UJLA (from ``unification'', ``Jordan'', 
``Lie'' and ``associative'').
structures, and
one could 
 ``decode'' the results obtained for UJLA structures
in results for  Jordan algebras, Lie algebras
or (non-unital) associative algebras.

Changhing the perspective, one can consider the  UJLA structures
as generalizations of Jordan algebras.

\begin{picture}(100,100)(10,10)
\put(60,80){$ \bf UJLA \ str. $ }
\put(70,10){$ \bf  \ k-alg $ }
\put(111,21){\vector(0,1){53}}
\put(240,80){$ \bf {UJLA \ str} $ }
\put(230,10){$ \bf { \ Jordan \ alg}$ }
\put(250,21){\vector(0,1){53}}
\put(128,86){\vector(1,0){104}}
\put(137,16){\vector(1,0){87}}
\put(174,92){$ \{ \ ,  \} $}
\put(95,42){ $ I $ }
\put(248,42){ $ I $ }
\put(174,22){$ \{ \ , \}  $}
\end{picture}

\bigskip

UJLA structures can also be interpretated as in an intermediat step in the
process of associationg a Lie algebra to an associative algebra (see
the picture below).

\bigskip




\begin{picture}(100,100)(10,10)
\put(60,80){$ \bf UJLA \ str. $ }
\put(70,10){$ \bf  \ k-alg $ }
\put(111,21){\vector(0,1){53}}
\put(240,80){$ \bf {UJLA \ str} $ }
\put(230,10){$ \bf { \ k-Lie \ alg}$ }
\put(250,74){\vector(0,-1){53}}
\put(128,86){\vector(1,0){104}}
\put(137,16){\vector(1,0){87}}
\put(174,92){$ [ \ ,  ] $}
\put(95,42){ $ I $ }
\put(248,42){ $ [ \ , ] $ }
\put(174,22){$ [ \ , ]  $}
\end{picture}

\bigskip

The study of filiform Lie algebras (\cite{ElisR})
can be extended to filiform  UJLA structures.

\bigskip

\begin{de}
We have defined the unifying structure $ ( V, \ \eta )$, 
also called a ``UJLA structure'', in the following way.
Let V be a vector space, and
   $ \eta : V \otimes V \rightarrow V, \ \ 
\eta (a \otimes b) = ab ,  $  be a linear map, which satisfies 
 the following axioms $ \forall a,b,c \in V$:
 
\begin{equation} \label{new}
 (ab)c + (bc)a + (ca)b  = a(bc) + b(ca) + c(ab),
\end{equation} 
\begin{equation} \label{Jordan1}
 (a^2 b) a \ = \ a^2 (ba), \ 
\end{equation} 
\begin{equation} \label{Jordan2}
(a b) a^2 \ = \ a (b a^2), \
\end{equation} 
\begin{equation} \label{Jordan3}
 (b a^2) a \ = \  (ba) a^2 , \
\end{equation} 
\begin{equation} \label{Jordan4}
 a^2 (ab) \ = \  a (a^2 b).
\end{equation} 

If just the identity (\ref{new}) holds, we call the
structure $ ( V, \ \eta )$ a ``weak unifying structure''.
\end{de}

\begin{re} If $ (A, \ \theta )$, where $ \theta : A \otimes A \rightarrow A, \ \ 
\theta (a \otimes b) = ab $, 
 is a (non-unital) associative algebra, then we define  a UJLA structure $ (A, \ \theta' )$, where $ \theta' (a \ot b) = \alpha ab \ + \ \beta ba $, for some
 $  \ \alpha, \ \beta \in k$.
For
$\alpha = \frac{1}{2}$ and $ \beta = \frac{1}{2}$, then $ (A, \ \theta' )$ is a Jordan algebra, and
for $\alpha = 1$ and $ \beta = -1$, then $ (A, \ \theta' )$ is a Lie algebra. 
\end{re}

\begin{te}
Let $(V, \eta) $ be a UJLA structure, and $ \ \alpha, \ \beta \in k$. Then,
 $(V, \eta'), \ \ \eta'(a\ot b)= \alpha ab+ \beta ba $ is a
 UJLA structure.

\end{te}

\begin{re}
Let $(V, \eta) $ be a UJLA structure. Then, $ \delta_a : V \rightarrow V, \ \  \delta_a (x) = xa-ax $,
is a derivation for the following UJLA structure:
 $(V, \eta'), \ \ \eta'(a\ot b)= ab - ba $. 
\end{re}

\bigskip
{\bf OPEN PROBLEM.}
A UJLA structure is power associative.
We know that a UJLA structure is power associative for dimensions less or equal to 5.

\begin{re}
 The classification of UJLA structures is also an open problem.
\end{re}

\begin{te}
Let $V$ be a vector space over the field $k$, and $p,q \in k$.
For $f,g : V \rightarrow V $, we define $ M (f \ot g) = f * g = f *_{p,q} g = p f \circ g + q g \circ f: V \rightarrow V$.
Then:

(i) $ ( End_k (V), \ *_{p,q}) $ is a UJLA structure $ \  \forall p,q \in k$.

(ii) For $ \phi : End_k (V) \rightarrow End_k (V \ot V)$ a morphism of UJLA structures (i.e., $ \ \phi (f*g)= \phi(f) * \phi(g) \ $),
$ \ W= \{ f: V \rightarrow V \vert f \circ M = M \circ \phi (f) \}$ is a sub-UJLA structure of the structure defined at (i). In other words,
$ \ f*g \in W, \ \forall f, g \in W $.

\end{te}

\begin{te}
Let $(V, \eta) $ be a UJLA structure. Then,
 $(V, \eta'), \ \ \eta'(a\ot b)=  ab - ba $ is a Lie algebra.
\end{te}
\begin{Proof}
 The proof follows from formula (\ref{new}).
\end{Proof}


\begin{te}
Let $(V, \eta) $ be a UJLA structure. Then,
 $(V, \eta'), \ \ \eta'(a\ot b)=  {\frac{1}{2}} (ab + ba) $ is a Jordan algebra.
\end{te}
\begin{Proof}
 The proof follows from formulas (\ref{Jordan1}), (\ref{Jordan2}), (\ref{Jordan3}) and (\ref{Jordan4}).
\end{Proof}


\begin{te}\label{tn}
For $(V, \eta) $ a UJLA structure, $ D (x) = D_b (x) = bx - xb$ 
is a { UJLA--derivation} (i.e., 
 $ D(a^2 a) =  D(a^2 )a + a^2 D(a) \ \forall a \in V$).
\end{te}

{\bf Proof.} In formula (\ref{new}) we take $ c= a^2$:
$ \ \ \  (ab)a^2 + (b  a^2 )a + (a^2 a)b  = a(ba^2 ) + b(a^2 a) + a^2 (ab) $.
It follows that
 $  (b  a^2 )a + (a^2 a)b  =  b(a^2 a) + a^2 (ab) $.

So, $  (b  a^2 )a  -  a^2 (ab) =  b(a^2 a) -  (a^2 a)b   $; so, 
$   b(a^2 a) -  (a^2 a)b = (b  a^2 - a^2 b )a  +  a^2 (ba - ab)   $.

Thus, $ D(a^2 a) =  D(a^2 )a + a^2 D(a)$.

\begin{de}
For the vector space V, let $ d: V \rightarrow V $ and
   $ \phi : V \ot V \rightarrow V \ot V,$ \ \ 
  be a linear map which satisfies:
\begin{equation} \label{ybptd}
 \phi^{12} \circ \phi^{23} \circ \phi^{12} =
 \phi^{23} \circ \phi^{12} \circ \phi^{23} 
\end{equation} 
where $ \ \phi^{12} = \phi \ot I, \ \phi^{23} = I\ot \phi $, $ \ \ I: V  \rightarrow V ,
 \ \ a  \mapsto a $.

Then, $(V, d, \phi) $ is called a { generalized derivation} if
$ \phi \circ ( d \ot I + I \ot d) = ( d \ot I + I \ot d) \circ \phi $.

\end{de}

\begin{re}
 If $A$ is an associtive algebra, $ d: A \rightarrow A $ a derivation
(so, $d(1_A)$=0),
and $ \phi : A \ot A \rightarrow A \ot A, \ a \ot b \mapsto ab \ot 1 +
1 \ot ab - a \ot b$, then $(A, d, \phi) $ is  a generalized derivation.

 If $C$ is a coalgebra, $ d: C \rightarrow C $ a coderivation,
and $ \psi : C \ot C \rightarrow C \ot C, \ c \ot d \mapsto 
\eps (d) c_1 \ot c_2 + \eps (c) d_1 \ot d_2
- c \ot d$, then $(C, d, \psi) $ is  a generalized derivation.

If $\tau$ is the twist map, the condition
$$ \tau \circ R \circ \tau = R $$ 
represents the unification of the
comutativity and the co-comutativity conditions.
In other words, if the algebra $A$ is comutative, then $ \phi$ verifies
the above condition. If the coalgebra $C$ is cocomutative, then $ \psi$ verifies
the same condition.
\end{re}


\begin{de}
 Let $A$ is an associtive algebra, $ d: A \rightarrow A $ a derivation,
$M$ an A-bimodule, and $ D:M \rightarrow M $ with the property
$ D(am) = d(a)m + a D(m)$.
Then, $( A, d, M, D)$ is called a { module derivation}.
\end{de}

\begin{te} ( \cite{uni2}) \label{der}
 In the above case, $A \times M$ becomes an algebra, and
$ \delta: A \times M \rightarrow A \times M, \ (a, \ m) \mapsto (d(a), D(m))$
is a derivation in this algebra.
\end{te}

Translated into the ``language'' of Differential Geometry, the above theorem 
says that
the Lie derivative is a derivation (i.e., $ \ d(ab)=d(a)b+ ad(b) \ $) on the product of 
 the algebra of functions defined on the manifold M with
 the set of vector fields on M (see \cite{geo}).

\begin{re} \label{coder}
 A dual construction would refer to a coalgebra structure,
$ \Delta : A \rightarrow A \otimes A, \ f \mapsto f \otimes 1 + 1 \otimes f$,
and a comodule structure on forms, 
$ \rho : \Omega \rightarrow A \otimes \Omega,
 \  f dx_1 \wedge dx_2 \ ... \wedge dx_n 
\mapsto f \otimes dx_1 \wedge dx_2 \ ... \wedge dx_n + 1 \otimes 
f dx_1 \wedge dx_2 \ ... \wedge dx_n $ .

In this case $A \oplus \Omega$ becomes a coalgebra.
\end{re}

\begin{re}
 The unification of Theorem \ref{der} and the construction
from Remark \ref{coder} can be realized at the level of 
Yang-Baxter systems.

Take, for example $ W(f \otimes g) = 1 \otimes fg$ and $X(f \otimes v)=
1 \otimes fv$ to obtain a semi Yang-Baxter systems.

\end{re}

{\centering\section{Conclusions}}

The Yang--Baxter systems were needed for unifying two constructions.

 Vr\u anceanu - Vergne basis were topics in a talk at DGA 13
 (see also \cite{ElisR}).

\bigskip

\section*{Acknowledgment}

We would like to thank Prof. Lazlo Stacho and Prof. Nicolae Anghel 
for their help and suggestions.
Also, we thank   the
Simion Stoilow Institute of Mathematics
of the Romanian Academy and the Petroleum Gas University from Ploiesti.

\end{document}